\documentclass{amsart}
\usepackage{amssymb,amsmath}
\usepackage[american,french]{babel}
\usepackage{amsopn}

\usepackage{fancyhdr}

\def\Titre{Spatial smoothness of the stationary solutions of the 3D Navier--Stokes equations}
\title{\Titre}
\author{Cyril ODASSO
\\
 \\  Ecole Normale Sup\'erieure de Cachan, antenne de Bretagne,\\ Avenue Robert Schuman,
 Campus de Ker Lann, 35170 Bruz (FRANCE). \\ and
\\ IRMAR,  UMR 6625 du CNRS, Campus de Beaulieu,  35042 Rennes cedex (FRANCE)\\}

\newtheorem{Theorem}{Theorem}[section]

\newtheorem{Lemma}[Theorem]{Lemma}
\newtheorem{Corollary}[Theorem]{Corollary}

\newtheorem{Hypothesis}[Theorem]{Hypothesis}
\makeatletter
\newif\ifmsbmloaded@

\@addtoreset{equation}{section}

\makeatother

\def\R{\mathbb R}
\def\N{\mathbb N}
\def\Z{\mathbb Z}

\def\E{\mathbb E}
\def\P{\mathbb P}

\def\H{\mathbb H}

\def\H{\Bbb H}

\newcommand{\BLANC}[1]{   }

\newcommand{\abs}[1]{\left\vert#1\right\vert}
\newcommand{\norm}[1]{\left\Vert#1\right\Vert}
\newcommand{\eps}{\varepsilon}

\def \Espace{\renewcommand{\arraystretch}{1.7} }

\begin{document}

\selectlanguage{american}

\maketitle

\pagestyle{fancy}

\noindent\textbf{Abstract}:
We consider stationary solutions of the three dimensional Navier--Stokes equations (NS3D)
with periodic boundary conditions and driven by
an external force which might have a deterministic and a random part.
The random part of the force is white in time and very smooth in space.
 We  investigate smoothness properties
in space of the stationary solutions.

Classical technics for studying smoothness of stochastic PDEs do not seem to apply since
global existence of strong solutions is not known. We use the
Kolmogorov operator and Galerkin approximations. We first assume that the  noise has spatial regularity
of order $p$ in the $L^2$ based Sobolev spaces, in other words that its paths are in $H^p$.
Then we prove that at each fixed time the law of the
stationary solutions is supported by $H^{p+1}$.

Then, using a totally different technic, we prove that if the noise has Gevrey regularity then
at each fixed time, the law of a stationary solution is supported by a Gevrey space.
Some informations on the Kolmogorov dissipation scale are deduced.


\section*{Introduction}

We are concerned with the stochastic Navier--Stokes equations in dimension $3$ (NS3D) with periodic boundary
 conditions and zero mean value.
These equations describe the time evolution of an incompressible fluid and are given by
\begin{equation}\label{Smooth_EqIntroNS}
\Espace
\left\{
\begin{array}{l}
dX+\nu(-\Delta)X \, dt+(X,\nabla)X\,dt+\nabla p \, dt =  \phi(X)dW+ g(X)dt,\\
\begin{array}{rcll}
        \left(\textrm{div } X\right)(t,\xi) &=& 0, &\textrm{ for } \xi\in D,\;\; t>0,\\
        \int_DX(t,\xi)d\xi &=& 0,&\textrm{ for } t>0,\\
        X(0,\xi)&=& x_0(\xi), &\textrm{ for } \xi\in D,
\end{array}
\end{array}
\right.
\end{equation}
where $D=(0,2\pi)^3$. 
We have denoted by $X(t,\xi)$ the velocity and by $p(t,\xi)$ the pressure at time $t$ and at the point
   $\xi\in D$, also $\nu$ denotes the viscosity.
 The  external force  acting on the fluid  is the sum of a
random force of white noise type $\phi(X)dW$ and a deterministic
one $g(X)dt$.

As it is well known, in the deterministic case, global existence
of weak (in the PDE sense) solutions and uniqueness of strong
solutions hold for the Navier-Stokes equations. In space dimension
two, weak solutions are strong and global existence and uniqueness
follows. Such a result is an open problem in dimension three  (see
\cite{temam-survey} for a survey on these questions).

In the stochastic case, the situation is similar. However due to
the lack of uniqueness, we have to work with global weak (in the
PDE sense) solutions of the martingale problem (see
\cite{flandoli-cetraro} for a survey on the stochastic case).
 Roughly speaking, this means
that in \eqref {Smooth_EqIntroNS}, we take $X$, $p$ and $W$ for unknown.

As is usual in the context of the incompressible Navier-Stokes equation, we get rid of the
pressure thanks to the Leray projector.
Let us denote by $(X,W)$ a weak (in the PDE sense) stationary
solution of the martingale problem
\eqref{Smooth_EqIntroNS} and by $\mu$ the 
 the law of $X(t)$, which is an invariant measure if we can prove that
 \eqref {Smooth_EqIntroNS} defines a Markov evolution. 
In this article, we establish that $\mu$ admits a finite moment in
spaces of smooth functions provided the external force  is
sufficiently smooth. We think that this is an interesting question
to study. First, it can be seen that if we were able to prove that
$\mu$ has a moment of sufficiently high order in a well chosen
Sobolev norm (order $4$ in $H^1$ or $2$ in $H^2$ for instance)
then this would imply global existence of strong solutions for
$\mu$ almost every initial data.

Moreover,  this result is an important ingredient if one tries to
follow the method of \cite{DebusscheNS3D} to construct a Markov
transition semi-group in $H^p(D)$  under suitable conditions on
$\phi$ and $g$. Since even uniqueness in law is not known for
NS3D, such a result might be important.

We  first prove that if the external force  is in $H^{p-1}(D)$ and the noise term
has paths in  $H^{p}(D)$ then
 $\mu$ admits a finite moment in the Sobolev space $H^{p+1}(D)$

Note that analogous results are well-known for the two dimensional Navier--Stokes equations (NS2D). Actually a stronger result
 is true for NS2D.
Namely, for any square integrable $x_0$, the unique solution of NS2D is
continuous from $(0,\infty)$ into $H^p(D)$ and is square
integrable from $(t_0,t_1)$ into
 $H^{p+1}(D)$. It follows that $\mu$ admits  moments of any orders in $H^{p}(D)$ and a moment of order 2 in $H^{p+1}(D)$.
This stronger result is linked to the global existence of  strong solutions for NS2D.

This kind of idea cannot be used for NS3D and we use  a generalization of an idea used
in \cite{DebusscheNS3D} for the case $p=1$. The method is based on the use
of the Kolmogorov operator applied to suitable Lyapunov functional. These functionals
have already been used in the deterministic case in \cite{Temam2}, chapter 4.

Using a totally different method, we establish also  that the invariant measure $\mu$ admits a moment
in a Gevrey class of functions provided the external force
 has the same regularity.
 Gevrey regularity has been studied in the deterministic case in  \cite{FoiasTemam} and \cite{Henshaw}.
 Our method is based on tools developed in \cite{FoiasTemam}. In  \cite{MGev}, \cite{SGevrey}  these tools
 have been used to obtain an exponential moment for the invariant measure
 in Gevrey norms in the two dimensional case. The arguments used in  \cite{MGev}, \cite{SGevrey}  do not generalize
 to the three dimensional case since there strong existence and uniqueness is used.
The three dimensional case NS3D requires
substantial adaptations.  We develop a framework which gives a control on a Gevrey norm
 by using a control of the $H^1(D)$--norm of $X$  only at fixed time.

 Actually, in this way,
we are able to generalize to NS3D the results of  \cite{MGev},
\cite{SGevrey}. However, we do not have exponential moments. We
deduce that the Kolmogorov dissipation scale is larger than
$\nu^{6+\delta}$. This is certainly not optimal since it is
expected that the scale is of order $\nu^{\frac 34}$. Note that
our result is rigorous and does not use any heuristic argument.

\section{Notations}

\noindent For $m\inÊ\N$, we denote by $\H^m_{\textrm{per}}( D)$
the space of functions which are restrictions of periodic
functions in $H^m_{loc}(D)^3$ and whose average
 is zero on $D$.
We set
$$
H=\left\{X\in \H^0_{\textrm{per}}(D)\,\left|\,\textrm{ div }X=0 \textrm{ on } D \right.\right\},
$$
and 
$$
V=H\cap \H^1_{\textrm{per}}(D).
$$

\noindent Let $\pi$ be the orthogonal projection in $L^2(D)^3$
onto the space $H$. We set
$$
A = \pi\left(-\Delta\right), \quad  D(A) =  V  \cap \H^2_{\textrm{per}}(D) \;\textrm{ and }\;
B(u)=\pi \left((u,\nabla)u\right).
$$
It is convenient to endow  $\H^m_{\textrm{per}}( D)$ with the
inner product $((\cdot,\cdot))_m=(A^{\frac m2}\cdot,A^{\frac
m2}\cdot)_{L^2(D)^3}$. The corresponding norm is denoted by
$\norm{\cdot}_m$. It is classical that this defines a norm which
is equivalent to the usual one. For $m=0$ we write
$|\cdot|=\norm{\cdot}_0$ and for $m=1$ we write
$\norm{\cdot}=\norm{\cdot}_1$. Note that, since we work with
functions whose average is zero on $(0,2\pi)^3$, we have the
following Poincar\'e type inequality
$$
\norm{x}_{m_1}\le  \norm{x}_{m_2},\; m_1\le m_2,\; x\in
\H_{\textrm{per}}^{m_2}(D).
$$
We also use the spaces $L^p(D)^3$ endowed with their usual norm
denoted  by $\abs\cdot_p$. Moreover, given two Hilbert spaces
$K_1$ and $K_2$,  $\mathcal{L}_2(K_1;K_2)$ is the space of
Hilbert-Schmidt operators from $K_1$ to $K_2$ .

The noise is described by a cylindrical Wiener process $W$ defined on a Hilbert space $U$ and
a mapping $\phi$ defined on  $H$ with values in $\mathcal L_2(U;H)$. We also consider
a deterministic forcing term described by a mapping $g$ from $H$ into $H$.
More precise assumptions on $\phi$ and $g$ are made below.

\noindent  Now, we can write problem \eqref{Smooth_EqIntroNS} in the form
\begin{equation}\label{Smooth_EqNS}
\Espace
\left\{
\begin{array}{rcl}
dX+  \nu A X dt+B(X) dt &=& \phi(X)dW+g(X)dt,\\
                             X(0)&=& x_0.
\end{array}
\right.
\end{equation}

\noindent In all the paper, we consider a $H$--valued stationary
solution $(X,W)$ of the martingale problem \eqref{Smooth_EqNS}.
Existence of such a solution has been proved in
\cite{FlandoliGatarek}.
 We denote by $\mu$
the law of $X(t)$. We do not consider any stationary solutions but
only those which are limit in  distribution of stationary
solutions of Galerkin approximations of \eqref{Smooth_EqNS}. More
precisely, for any $N\in\N$, we denote by $P_N$ the projection of
$A$ onto the vector space spanned by the first $N$ eigenvalues and
consider the following approximation of \eqref{Smooth_EqNS}
\begin{equation}\label{Galerkin}
\Espace
\left\{
\begin{array}{rcl}
dX_N+  \nu A X_N dt+P_NB(X_N) dt &=& P_N\phi(X_N)dW+P_Ng(X_N)dt,\\
                             X_N(0)&=& P_N x_0.
\end{array}
\right.
\end{equation}
It can be easily shown that \eqref{Galerkin} has a stationary
solution $X_N$. Proceeding as in \cite{FlandoliRomito}, we can see
that their laws are tight in suitable functional spaces, and, up
to a subsequence, $(X_N,W)$ converges in law to a stationary
solution $(X,W)$ of \eqref{Smooth_EqNS}. Actually the convergence
holds in $C(0,t;D(A^{-s}))\cap L^2(0,t;D(A^{\frac 12-s}))$ for any
$t,s>0$. We only consider stationary solutions constructed in that
way.

To obtain an estimate for stationary solutions of
\eqref{Smooth_EqNS} (limit of Galerkin approximations), we proceed
as follows. We first prove the desired estimate for every
stationary solutions of \eqref{Galerkin} and then we take the
limit.

The reason why our results are only stated for solutions limit of
Galerkin approximations comes from the fact that it is not known
if computations applied to solutions of Galerkin approximations
can be applied directly on solutions $(X,W)$ of the
three-dimensional Navier-Stokes equations.

\noindent Some of our results describe properties of $\mu$ in
Gevrey type spaces. These spaces contain functions with
exponentially decaying Fourier coefficients. According to the
setting given in \cite{FoiasTemam},
 we set for any $(\alpha,\beta)\in \R^+_*\times(0,1]$
$$
\Espace
\left\{
\begin{array}{rcl}
 \norm{x}_{G(\alpha,\beta)}^2&=&\abs{A^\frac{1}{2}e^{\alpha A^\frac{\beta}{2}}x}^2
 =\sum_{k\in Z^3}\abs{k}^2e^{2\alpha
\abs{k}^\beta}\abs{\hat x(k)}^2,\\
G(\alpha,\beta)&=&\left\{ x \in H\,\left|\,
\norm{x}_{G(\alpha,\beta)}<\infty  \right.\right\},
\end{array}
\right.
$$
where $(\hat x(k))_{k\in \Z^3}$ are the Fourier coefficients of
$x\in H$. Moreover, for any $(x,y)\in G(\alpha,\beta)^2$, we set
$$
\left(x,y\right)_{G(\alpha,\beta)}= \left(A^\frac{1}{2}e^{\alpha
A^\frac{\beta}{2}}x,A^\frac{1}{2}e^{\alpha
A^\frac{\beta}{2}}y\right)=\sum_{k\in Z^3}\abs{k}^2e^{2\alpha
\abs{k}^{\beta}}\Re e\left(\hat x(k)\overline{\hat y(k)}\right).
$$
 Clearly,
$\left(G(\alpha,\beta),\left(\cdot,\cdot\right)_{G(\alpha,\beta)}\right)$
is a
  Hilbert space.

We are not interested in large viscosities and in all the article it is assumed that $\nu\leq 1$.
We will use various constants which may depend on some parameter such as $p,\nu,\dots$ When this dependance
 is important, we make it explicit.

\section{$\H^p_{\textrm{per}}(D)$--regularity}


\noindent Let $p\in\N$. We now make the following smoothness assumptions on the forcing terms.

\
\begin{Hypothesis}
\label{Hp} The mapping $\phi$ (resp. $g$) takes values in
$\mathcal L_2\left(U;H\cap \H^p_{\textrm{per}}(D) \right)$ (resp.
$H\cap \H^{p-1}_{\textrm{per}}(D)$) and $\phi:H\to \mathcal
L_2\left(U;H\cap \H^p_{\textrm{per}}(D) \right)$ and $g:H\to H\cap
\H^{p-1}_{\textrm{per}}(D)$ are bounded.
\end{Hypothesis}


\

\noindent We set, when Hypothesis \ref{Hp} holds,
$$
B_p=\sup_{
H}\left(\norm{\phi}^2_{\mathcal{L}_2(U;\H^p_{\textrm{per}}(D))}+\norm{g}_{p-1}^2\right).
$$
It is also convenient to define
$$
\bar B_p=\sup_{
H}\norm{\phi}^2_{\mathcal{L}_2(U;\H^p_{\textrm{per}}(D))}+
\frac{1}{\nu}\sup_{ H}\norm{g}_{p-1}^2.
$$
The aim of this section is to establish the following result.
\begin{Theorem}\label{Smooth_ThMAIN}
Let $\mu$ be the invariant law of a stationary solution $X$ of the
three dimensional Navier-Stokes equations that is limit of
stationary solutions of Galerkin approximations. Assume that
Hypothesis \ref{Hp} holds for some $p\ge 1$.  For any $\nu\leq 1$,
there exists $c_{p,\nu}$ depending on $p$, $\nu$ and $B_p$ such
that
$$
\int_H \norm{x}_{p+1}^{\frac{2}{2p+1}}\,d\mu(x)\leq c_{p,\nu}.
$$
\end{Theorem}
Let us make few comments.

Note that it would be very important to obtain an estimate on $\int_H \norm{x}_{p+1}^{\delta_p}\,d\mu(x)$
with $ p\delta_p>3$. Indeed, by Agmon  inequality , we have
$$
\int_H |x|_{\infty}^{2}\,d\mu(x)\le c \int_H |x|^{2-3/p} \norm{x}_{p}^{\frac{3}{p}}\,d\mu(x)
$$
and this would give an estimate on the left hand side. Since
uniqueness is easily shown to hold for solutions in $L^2(0,T;
L^\infty(D)^3)$, a classical argument could be used to deduce that
for $\mu$ almost every initial data there exists a unique global
weak solution. Combining with the result in
\cite{flandoli-irreducibility}, this would partially solve Leray's
conjecture.

Consider the case $g=0$, $U=H$ and $\phi=A^{-s-\frac{3}{4}}$. Then
Hypothesis \ref{Hp} holds for any $p<s$ and the
 unique invariant measure of the three dimensional linear stochastic Stokes equations
  in $H$ is in $\H^{r+1}_{\textrm{per}}(D)$
 with probability zero if $r>s$.
  Therefore it seems that $\norm\cdot_{p+1}$ is
the strongest norm we can control under Hypothesis \ref{Hp}.

Remark that in the two dimensional case a much stronger result holds. Indeed, standard arguments
imply that under Hypothesis \ref{Hp} we have
 for any invariant measure $\mu$ and any $q\in\N^*$
$$
\int_H \norm{x}_{p}^{2q}\,d\mu(x)<\infty,\quad \int_H \norm{x}_{p+1}^{2}\,d\mu(x)<\infty.
$$

In the proof, we use ideas developped in \cite{Temam2}. Similar
but more refined techniques have been used in \cite{Henshaw} to
derive interesting properties on the decay of the Fourier spectrum
of smooth solutions of the deterministic Navier-Stokes equations.
Using such techniques does not seem to yield great improvement of
our result. Indeed, trying to do so, we have been able to improve
the estimate of Theorem \ref{Smooth_ThMAIN} as follows
$$
 \nu\int_H \norm{x}_{p+1}^{\frac{c_*}{p}}\,d\mu(x)\leq 2\bar B_p+ c 2^p (1+\bar B_0),
$$
 where $c$ and $c_*$ are positive constants and $c_*$ is close to $1.02$.
We have not been able to derive very interesting results from this improved estimate and therefore
have preferred to give the simpler one which follows from easier arguments.




{\bf Proof}: The proof is rather standard. We do not give the
details.

 Let $(\mu_{N})_{N\in\N}$ be a sequence of invariant
measures of stationary solutions $(X_N)_N$ of \eqref{Galerkin}
such that there exists a subsequence $(N_k)_{k\in \N}$ such that
$X_{N_k}$ converges to $X$ in law. It follows that
$(\mu_{N_k})_{k\in\N}$  converges to $\mu$ (considered as
probability measures on $D(A^{-1})$).

\noindent We denote by $L_N$ the  Kolmogorov operator associated
to the Galerkin approximation \eqref{Galerkin} of the stochastic
Navier-Stokes equations
$$
L_Nf(x)=\frac{1}{2}tr\left((P_N\phi)(x)(P_N\phi)^*(x)D^2f(x)\right)-\left(\nu
Ax+B(x)-g(x),Df(x)\right),
$$
for any $f\in C^2(P_N H; \R)$ and $x\in P_N H$.

\noindent The proof of Theorem \ref{Smooth_ThMAIN} is based on the
fact that,  for any $N\in \N$, we have
\begin{equation}\label{Smooth_EqKol1}
\int_{P_NH} L_Nf(x) d\mu_N(x)=0,
\end{equation}
provided $f\in C^2(P_N H; \R)$  verifies
\begin{equation}\label{Cond1}
\Espace \left\{
\begin{array}{rrcl}
i)& \int_{P_NH} \abs{f(x)}d\mu_N(x)&<&\infty,\\
ii)& \int_{P_NH} \abs{L_Nf(x)}d\mu_N(x)&<&\infty,\\
iii)& \int_{P_NH} \abs{(P_N\phi)^*(x)Df(x)}^2d\mu_N(x)&<&\infty.
\end{array}\right.
\end{equation}

\BLANC{{\bf Step 0:} A preliminary result.

It\^o formula of $\abs{X_N}^{2}$ gives
$$
d\abs{X_N}^{2}+2\nu\norm{X_N}^2dt=2(X_N,\phi(X_N)dW)+2(X_N,g(X_N))dt+\norm{P_N\phi(X_N)}^2_{\mathcal
L_2(U;H)}dt.
$$
Taking into account $\abs\cdot\leq \norm\cdot$ and
$2(X_N,g(X_N))\leq \nu \norm{X_N}^2+\frac 1\nu
\norm{g(X_N)}_{-1}^2dt$, it follows that
$$
d\abs{X_N}^2+\nu\abs{X_N}^2dt\leq 2(X_N,\phi(X_N)dW)+\bar B_0 dt.
$$
It\^o formula of $\abs{X_N}^{2p}$ gives
$$
\Espace \begin{array}{r}
 d\abs{X_N}^{2p}+\nu p\abs{X_N}^{2p}dt\leq
2p\abs{X_N}^{2(p-1)}(X_N,\phi(X_N)dW)+p\bar B_0 \abs{X_N}^{2(p-1)}
dt\quad\quad\\
+\frac{p(p-1)}{2}\abs{\phi^*(X_N)X_N}^2\abs{X_N}^{2(p-2)}dt.
\end{array}
$$
Arithmetic-geometric inequality gives
$$
d\abs{X_N}^{2p}+\nu\abs{X_N}^{2p}dt\leq
2p\abs{X_N}^{2(p-1)}(X_N,\phi(X_N)dW)+C_{p,\nu}dt.
$$
 Integrating on $(-t,0)$ and
taking the expectation, we obtain for any $x\in P_N H$ and $t>0$
$$
\E\left(\abs{X_N(0)}^{2p}\,\left|\,X_N(-t)=x\right.\right)\leq
e^{-\nu t}\abs{x}^{2p}+C_{p,\nu},
$$
which yields
\begin{equation}\label{Eq0.7}
\E\left(\abs{X_N(0)}^{2p}1_{\abs{X_N(-t)}^{2p}\leq M}\right)\leq
e^{-\nu t}M+C_{p,\nu}.
\end{equation}
 For any $t>0$, we set $M_t=e^{\nu t}C_{p,\nu}$. The law of $X_N(-t)$
 is not depending on $t$. Hence, taking the limit, it follows, by
  classical arguments, that there exists a sequence $(t_k)_{k\in\N}$
   tending to infinity such that
   $$
1_{\abs{X_N(-t_k)}^{2p}> M_{t_k}}\to 0\,\textrm{ a.s.}
   $$
   Taking the limit in \eqref{Eq0.7}, it follows  from
Fatou Lemma}

It follows from \cite{flandoli-cetraro}, Chapter 1.2, Corollary
1.12 that for any $p\in\N$
\begin{equation}\label{Eq0}
\int_{P_N H}\abs{x}^{2p}d\mu_N(x)=
\E\left(\abs{X_N(0)}^{2p}\right)\leq 2C_{p,\nu}<\infty,
\end{equation}
and
\begin{eqnarray}
\nu \int_{P_NH} \norm x ^2 d\mu_N(x)&\leq&\bar B_0 .\label{Smooth_Eq2.1.1}
\end{eqnarray}
The result in \cite{flandoli-cetraro} is given for $g=0$ but the
generalization is easy.

 Thanks to \eqref{Eq0}, we use
\eqref{Smooth_EqKol1} with
$$
f=\frac{1}{\left(1+\norm \cdot ^2_{p}\right)^{\eps_p }},\quad
\eps_p= \frac{1}{2p-1}.
$$
We obtain after lengthy but easy computations and some estimates
on the nonlinear term borrowed from \cite{Temam2}, chapter 4 (in
particular equation (4.8)) that  there exists $c_p$ such that
\begin{equation}\label{Smooth_Eq2.2.0}
 R_p\leq 2\bar B_p+c_{p}\bar B_0 +1,
\end{equation}
with
$$
R_p=\nu\int_{P_NH} \frac{1+\norm x ^2_{p+1}\quad}{\left(1+\norm x
^2_{p}\right)^{1+\eps_p }} d \mu_N(x).
$$
Then arguing as in \cite{Temam2}, chapter 4, we set
$$
M_p=\nu \int_{P_NH} \left(
1+\norm{x}_{p}^2\right)^{1/2p-1}d\mu_N(x),
$$
 and deduce
$$
M_{p+1}\le R_p^{1/2p+1} M_p^{2p/2p+1},
$$
which yields
\begin{equation}\label{Eq1}
\int_{P_NH} \norm{x}_{p+1}^{\frac{2}{2p+1}}\,d\mu_N(x)\leq
c_{p,\nu}.
\end{equation}
It is then standard to deduce the results thanks to \eqref{Eq1},
the subsequence $(N_k)_k$, Fatou Lemma and lower semicontinuity of
$\norm\cdot_p$.

\section{Gevrey regularity}

\subsection{Statement of the result}

\noindent We now state and prove the main result of this paper. It
states that if
 the external force  is bounded in a Gevrey class of functions,
 then $\mu$ have support in another Gevrey class of function.

\noindent The main assumption in this section is the following

\begin{Hypothesis}
\label{G0} There exists $(\alpha,\beta)\in \R^*_+\times(0,1]$ such
that the mappings $g: H \to  G(\alpha,\beta)$ and   $\phi: H \to
\mathcal L_2\left(U;G(\alpha,\beta)\right)$ are bounded.
\end{Hypothesis}

\smallskip

\noindent We set
$$
B_0'=\sup_{x\in H}\norm{\phi(x)}^2_{ \mathcal L_2\left(U;G(\alpha,\beta)\right)}+
\sup_{x\in H}\norm{g(x)}^2_{G(\alpha,\beta)}.
$$

\noindent The aim of this section is to establish the following results proved in the
following subsections.

\begin{Theorem}\label{Gevrey_Th_Gevrey}
Let $\mu$ be the invariant law of a stationary solution $X$ of the
three dimensional Navier-Stokes equations that is limit of
stationary solutions of Galerkin approximations. Assume that
Hypothesis \ref{G0} holds. There exist a family of constants
$(K_\gamma)_{\gamma\in (0,1)}$ only depending on
$(\alpha,\beta,B_0')$ and a family $(\alpha_\nu)_{\nu\in (0,1)}$
of measurable mappings $ H\to(0,\alpha)$  such that for any
$\nu\in(0,1)$
\begin{eqnarray}
\int_{H} \norm{x}^{2\gamma}_{G( \nu \alpha_\nu(x),\beta)}d\mu(x)&\leq& K_\gamma(1+\bar B_0)^2\nu^{-\frac72}
,\label{Gevrey_Eq_Th_Gevrey_1}\\
\int_H \left( \alpha_\nu(x)\right)^{-\frac\gamma2}d\mu(x)&\leq& K_\gamma(1+\bar B_0)\nu^{-\frac52}
,\label{Gevrey_Eq_Th_Gevrey_2}
\end{eqnarray}
for any $\gamma\in(0,1)$.
\end{Theorem}

\noindent This result gives some informations on the Kolmogorov dissipation scale. Indeed, it
 follows from \eqref{Gevrey_Eq_Th_Gevrey_1}, \eqref{Gevrey_Eq_Th_Gevrey_2} that
$$
\abs{\hat x(k)}\leq \norm{x}_{G( \nu
\alpha_\nu(x),\beta)}\abs{k}^{-1}
e^{-\nu\alpha_\nu(x)\abs{k}^\beta},
$$
where $(\hat x(k))_{k\in\Z^3}$ are the Fourier coefficients of $x$.

\noindent Hence, if Hypothesis \ref{G0} holds with $\beta=1$ and
$g=0$, then $\abs{\hat x(k)}$ decreases faster than any powers
 of $\abs{k}$ for $\abs{k}>\hspace{-0.1cm}>(\nu\alpha_\nu(x))^{-1}$. By  \eqref{Gevrey_Eq_Th_Gevrey_2}, for any
$\delta>0$
$$
\frac{1}{\alpha_\nu(x)}\leq c_{\delta,\nu}(x)\nu^{-5(1+\delta)} \;
\textrm{ with } \; \int
c_{\delta,\nu}(x)^{\frac{1}{2(1+\delta)}}\,\mu(dx)\leq
\Theta_{\delta}<\infty,
$$
and $\Theta_\delta$ not depending on $\nu$. It follows
that $\abs{\hat x(k)}$ decreases faster than any powers
 of $\abs{k}$ for $\abs{k}>\hspace{-0.1cm}>\nu^{-(6+5\delta)}$.
This indicates that the 3D--Kolmogorov dissipation scale is larger
than $\nu^{6+5\delta}$. Note that by physical arguments it is
expected that the 3D--Kolmogorov dissipation scale is of order of
$\nu^{\frac 34}$.

Making analogous computations, we obtain that, for $g\not=0$ and
$\beta\in (0,1]$, the 3D--Kolmogorov dissipation scale is larger
than $\nu^{\frac 8\beta+\delta}$ for $\delta>0$.

In \cite{MGev}, \cite{SGevrey}, analogous results  to Theorem
\ref{Gevrey_Th_Gevrey}  are proved for the 2D--Navier-Stokes
equations. Moreover, it is deduced in \cite{MGev} that the
2D--Kolmogorov dissipation scale is larger than
$\nu^{\frac{25}{2}}$. In \cite{BKL00}, it is established that the
2D--Kolmogorov dissipation scale is larger than
$\nu^{\frac{1}{2}}$, which is the physically expected result.

\noindent We can also control a moment of a fixed Gevrey norm.
\begin{Corollary}\label{Gevrey_Cor_Gevrey}
Under the same assumptions, there exists a family $(C_{\gamma,\alpha',\beta',\nu})_{\gamma,\alpha',\beta',\nu}$
 only depending on $(\alpha,\beta,B_0',\nu)$ such that
\begin{equation}
\int \left(\ln^+ \norm{x}^2_{G(\alpha',\beta')}\right)^{\gamma}d\mu(x)\leq C_{\gamma,\alpha',\beta',\nu},\label{Gevrey_Eq_Cor_Gevrey}
\end{equation}
where $\ln^+r = \max\{0,\ln r\}$ and provided $\,\alpha'>0$ and $\beta',\gamma>0$ verify
$$
\beta'<\beta\quad\textrm{ and }\quad
2\gamma<\frac{\beta}{\beta'}-1.
$$
\end{Corollary}

\subsection{Estimate of blow-up time in Gevrey spaces}

\

\noindent Before proving Theorem \ref{Gevrey_Th_Gevrey}, we establish the following result 
 which implies that the time of blow-up of the solution in Gevrey spaces
 admits a negative moment. 
 \begin{Lemma}
 \label{Gevrey_lem_Gevrey}
Assume that Hypothesis \ref{G0} holds. For any $N$, any stationary
solution $X_N$ of \eqref{Galerkin} and any $\nu\in(0,1)$, there
exist $K$ only depending on $(\alpha,\beta,B_0')$ and a stopping
time
 $\tau^N>0$ such that
\begin{eqnarray}
\E\left( \sup_{t\in(0,\tau^N)}\norm{X_N(t)}^2_{G(\nu   t,\beta)}\right)&\leq& \frac 4\nu(\bar B_0+1),\label{Gevrey_Eq_lem_Gevrey_1}\\
\P\left(\tau^N< t\right)&\leq & K (\bar B_0+1)t^\frac{1}{2}\nu^{-\frac52}
.\label{Gevrey_Eq_lem_Gevrey_2}
\end{eqnarray}
\end{Lemma}

\noindent This result is a refinement of the result developed by
Foias and  Temam in \cite{FoiasTemam} and is strongly based on the
tools developed in this latter  article. We denote by $\mu_N$ the
invariant law associated to $X_N$. Let us set
\begin{equation}\label{Gevrey_Eq3.0}
\tau^N=\inf\left\{t\geq 0\,\left|\,
 1+\norm{X_N(t)}_{G(\nu   t,\beta)}^2>4\left(\norm{X_N(0)}^2+1\right)
   \right.\right\}.
\end{equation}
Clearly
$$
\E\left( \sup_{t\in(0,\tau^N)}\norm{X_N(t)}^2_{G(\nu
t,\beta)}\right)\le 4\E\left(\norm{X_N(0)}^2+1 \right)
$$
and
 \eqref{Gevrey_Eq_lem_Gevrey_1} follows from 
  \eqref{Smooth_Eq2.1.1}. It remains to prove \eqref{Gevrey_Eq_lem_Gevrey_2}.

\noindent We apply It\^ o formula to $\norm{X_N(t)}^2_{G(\nu
t,\beta)}$ for $t\in (0,\alpha)$
\begin{equation}\label{Gevrey_Eq3.2.1}
\Espace \begin{array}{r}
 d\norm{X_N(t)}^2_{G(\nu
t,\beta)}+2\nu\norm{A^{\frac{1}{2}}X_N(t)}^2_{G(\nu   t,\beta)}dt=
 \nu\norm{A^\frac{\beta}{2}X_N(t)}^2_{G(\nu   t,\beta)}dt
\quad\quad\\
+dM(t)+I(t)dt,
\end{array}
\end{equation}
 where
$$
\Espace
\left\{
\begin{array}{lcllcl}
I(t)            &=& 2I_g(t)+2I_B(t)+I_\phi(t),&
I_B(t)      &=& -\left(X_N(t),B(X_N(t))\right)_{G(\nu   t,\beta)},\\
I_\phi(t)   &=& \norm{P_N\phi(X_N(t))}_{\mathcal L_2\left(U;G(\nu
t,\beta)\right)}^2,&
I_g(t)      &=& \left(g(X_N(t)),X_N(t)\right)_{G(\nu   t,\beta)},\\
\lefteqn{M(t)       \;\;\,  =\,\;\; 2\int_0^t\left(X_N(s)
,\phi(X_N(s))dW(s)\right)_{G(\nu   t,\beta)}.}
\end{array}
\right.
$$
The following inequality is proved in \cite{FoiasTemam} for
$\beta\leq 1$
\begin{equation}\label{Gevrey_Eq3.2.4}
2 I_B(t)\leq \nu\norm{A^\frac{1}{2} X_N(t)}_{G(\nu
t,\beta)}^2+\frac{c}{\nu^3}\norm{X_N(t)}_{G(\nu   t,\beta)}^6.
\end{equation}
By Hypothesis \ref{G0} we have
\begin{equation}\label{Gevrey_Eq3.2.5}
I_\phi(t)+2I_g(t)\leq \norm{ X_N(t)}_{G(\nu   t,\beta)}^6+B_0'+1.
\end{equation}
Combining \eqref{Gevrey_Eq3.2.1}, \eqref{Gevrey_Eq3.2.4} and  \eqref{Gevrey_Eq3.2.5}, we obtain
since $\beta,\nu\le 1$
\begin{equation}\label{Gevrey_Eq3.2.6}
d\norm{X_N(t)}^2_{G(\nu   t,\beta)}
\leq dM(t)+\frac{c}{\nu^3}\norm{X_N(t)}_{G(\nu
t,\beta)}^6dt+(B_0'+1)dt.
\end{equation}
 Applying Ito formula to $\left(1+\norm{X_N(t)}^{2}_{G(\nu   t,\beta)}\right)^{-2}$, we then deduce from \eqref{Gevrey_Eq3.2.6} and
from Hypothesis \ref{G0} that for
any  $t\in (0,\alpha)$ and any $\nu\leq 1$
\begin{equation}\label{Gevrey_Eq3.2.7}
-d\left(1+\norm{X_N(t)}^{2}_{G(\nu   t,\beta)}\right)^{-2}
\leq
d\mathcal M(t)+C_0\nu^{-3} dt,
\end{equation}
where $C_0=c(B_0'+1)$ and
$$
\mathcal M(t)=4{\int_0^t}\left(1+\norm{X_N(s)}^2_{G(\nu
t,\beta)}\right)^{-3} \left(X_N(s)
\,,\,\phi(X_N(s))dW(s)\right)_{G(\nu   t,\beta)} .
$$
Setting
$$
\Espace
\left\{
\begin{array}{rcl}
\tau^N_0&=&\inf\left\{t\in (0,\alpha)\,\left|\,
 \mathcal M(t)>\frac{1}{4\left(1+\norm{X_N(0)}^2\right)^{2}}
   \right.\right\},\\
\tau^N_1&=&\tau^N_0\wedge\left(\frac{\nu^{3}}{4C_0\left(1+\norm{X_N(0)}^2\right)^{2}}\right)
,
\end{array}
\right.
$$
we obtain by integration of \eqref{Gevrey_Eq3.2.7} on $[0,t]$ for
$t\in (0,\tau^N_1)$
$$
1+\norm{X_N(t)}^2_{G(\nu t,\beta)}\leq
4\left(1+\norm{X_N(0)}^2\right).
$$
We deduce that $\tau^N\geq \tau^N_1$ and
\begin{equation}\label{Gevrey_Eq3.2.9}
\P\left(\tau^N< t\right)\leq \P\left( \tau^N_0< t\right)+
\P\left(\left(1+\norm{X_N(0)}^2\right)^{2}\geq \frac{\nu^{3}}{4C_0
t}\right).
\end{equation}
Since $\mu$ is the law of $X_N(0)$, we have
$$
\P\left(\left(1+\norm{X_N(0)}^2\right)^{2}\geq \frac{\nu^{3}}{4C_0
t}\right)= \mu_N\left(x\in H , \,1+\norm{x}^2\geq
\frac{\nu^{\frac{3}{2}}}{(4C_0 t)^{\frac{1}{2}}}\right).
$$
Applying  Chebyshev inequality, we deduce from
\eqref{Smooth_Eq2.1.1}
\begin{equation}\label{Gevrey_Eq3.2.10}
\begin{array}{l}
\displaystyle{\P\left(\left(1+\norm{X_N(0)}^2\right)^{2}
\geq \frac{\nu^{3}}{4C_0 t}\right)}\\
\displaystyle{\leq
2\nu^{-\frac{3}{2}}\left(C_0 t\right)^\frac12\int_H(1+\norm{x}^2)d\mu_N(x)}\\
\\
\displaystyle{\le 2\nu^{-\frac{5}{2}}\left(1+\bar
B_0\right)\left(C_0 t\right)^\frac1{2} .}
\end{array}
\end{equation}
Moreover
$$
\begin{array}{ll}
\P\left( \tau^N_0< t\right)&=
\P\left(4\left(\left(1+\norm{X_N(0)}^2\right)^{2}\sup_{s\in[0,t\wedge \tau_0^N]}\mathcal M(s)\right)> 1\right)\\
& \le
4\E\left(\left(1+\norm{X_N(0)}^2\right)^2\sup_{s\in[0,t\wedge
\tau_0^N]}\mathcal M(s)\right).
\end{array}
$$
 Taking conditional expectation with respect to the $\sigma$--algebra $\mathcal F_0$  generated by $X_N(0)$ inside
  the expectation, it follows
$$
\P\left( \tau^N_0< t\right)\leq
4\E\left(\left(1+\norm{X_N(0)}^2\right)^2
\E\left(\left.\sup_{s\in[0,t\wedge \tau_0^N]}\mathcal
M(s)\,\right|\,\mathcal F_0\right) \right).
$$
By Burkholder-Davis-Gundy inequality (see Theorem 3.28 page 166 in
\cite{KARATZAS}) we obtain
$$
\E\left(\left.\sup_{s\in[0,t\wedge \tau_0^N]}\mathcal
M(s)\,\right|\,\mathcal F_0\right) \leq c
\E\left(\left.\left<\mathcal M\right>^\frac{1}{2}(t\wedge
\tau_0^N)\,\right|\,\mathcal F_0\right),
$$
and
\begin{equation}\label{Gevrey_Eq3.2.11}
\P\left( \tau^N_0< t\right)\leq
4\E\left(\left(1+\norm{X_N(0)}^2\right)^2
\E\left(\left.\left<\mathcal M\right>^\frac{1}{2}(t\wedge
\tau_0^N)\,\right|\,\mathcal F_0\right) \right).
\end{equation}
We have
$$
\left<\mathcal M\right>(t)= 4\int_0^t
\left(1+\norm{X_N(s)}^2_{G(\nu   t,\beta)}\right)^{-6}
\abs{\left(A^\frac{1}{2}e^{  \nu t
A^\frac{\beta}{2}}\phi(X_N(s))\right)^* \left(A^\frac{1}{2}e^{ \nu
t A^\frac{\beta}{2}}X_N(s)\right)}_U^2 ds.
$$
Therefore
$$
\left<\mathcal M\right>(t) \leq 4\int_0^t
\left(1+\norm{X_N(s)}^2_{G(\nu   t,\beta)}\right)^{-6}
\norm{\phi(X_N(s))}^2_{\mathcal L\left(U;G(\nu
t,\beta)\right)}\norm{X_N(s)}_{G(\nu t,\beta)}^2 ds .
$$
It follows that
$$
\left<\mathcal M\right>(t\wedge \tau_0^N) \leq \frac{B_0' t}{4^3
\left(1+\norm{X_N(0)}^2\right)^{4} }
$$
 Hence we infer from \eqref{Gevrey_Eq3.2.11} and from $\norm\cdot\leq\norm\cdot_{G(\nu   t,\beta)}$ that
\begin{equation}\label{Gevrey_Eq3.2.12}
\P\left( \tau^N_0\leq t\right)\leq   \sqrt{B_0't}.
\end{equation}
Combining \eqref{Gevrey_Eq3.2.9}, \eqref{Gevrey_Eq3.2.10} and \eqref{Gevrey_Eq3.2.12},
 we deduce \eqref{Gevrey_Eq_lem_Gevrey_2}.

\subsection{Proof of Theorem \ref{Gevrey_Th_Gevrey}}

\

\noindent Let $(\mu_{N})_{N\in\N}$ be a sequence of invariant
measures of stationary solutions $(X_N)_N$ of \eqref{Galerkin}
such that there exists a subsequence $(N_k)_{k\in \N}$ such that
$X_{N_k}$ converges to $X$ in law. It follows that
$(\mu_{N_k})_{k\in\N}$  converges to $\mu$ (considered as
probability measures on $D(A^{-1})$).

\noindent Setting
$$
\alpha_\nu(x)=\inf\left\{ s\geq 0  \,\left|\,  \norm{x}_{G(\nu
s,\beta)}^2> \frac{4}{\nu s^\frac{1}{2}}\left(\bar B_0+1\right)
 \right.\right\},
$$
it follows that for any $\gamma\in(0,1)$
\begin{equation}\label{Gevrey_Eq2.1}
\int \norm{x}^{2\gamma}_{G(\nu   \alpha_\nu(x),\beta)}d\mu(x)\leq
\left(\frac{4}{\nu}\right)^\gamma\left(\bar B_0+1\right)^\gamma
\int \left(\alpha_\nu(x)\right)^{-\frac{\gamma}{2}}d\mu(x).
\end{equation}
Hence \eqref{Gevrey_Eq_Th_Gevrey_1} is consequence of
\eqref{Gevrey_Eq_Th_Gevrey_2}. Then in order to establish Theorem
 \ref{Gevrey_Th_Gevrey}, it is sufficient to prove  \eqref{Gevrey_Eq_Th_Gevrey_2}.

\noindent Clearly
$$
\Espace \begin{array}{r}
 \P\left(\norm{X_N( t)}_{G(\nu
t,\beta)}^2> \frac{4}{\nu t^\frac{1}{2}}\left(\bar B_0+1\right)
\right)\leq \P\left(\sup_{s\in[0,\tau^N]}\norm{X_N(s)}_{G(\nu
s,\beta)}^2>\frac{4}{\nu t^\frac{1}{2}}\left(\bar
B_0+1\right)\right)\\
+ \P\left(\tau^N<t\right),
 \end{array}
$$ where $\tau^N$ has been
defined in section $3.2$. Applying  Chebyshev inequality, we infer
from Lemma \ref{Gevrey_lem_Gevrey}
 and from the fact that, for any $t>0$, $\mu_N$ is the law of $X_N(t)$
\begin{equation}\label{Gevrey_Eq2.3bis}
\mu_N(\mathcal O_t)=\P\left(X_N(t)\in \mathcal O_t\right)\leq (K
+1)(1+\bar B_0)
t^\frac12\nu^{-\frac52}.
\end{equation}
where
$$
\mathcal O_t=\left\{x\in D(A^{-1})\, ,\,\norm{x}_{G(\nu
t,\beta)}^2> \frac{4}{\nu t^\frac{1}{2}}\left(\bar
B_0+1\right)\right\}.
$$
Notice that $\mathcal O$ is an open subset of $D(A^{-1})$. Hence,
since $\mu_{N_k}\to \mu$ (considered as probability measures on
$D(A^{-1})$), then we deduce from \eqref{Gevrey_Eq2.3bis}
 that
\begin{equation}\label{Gevrey_Eq2.3}
\mu(\mathcal O_t)\leq \lim\inf_k\left(\mu_{N_k}(\mathcal
O_t)\right) \leq (K +1)(1+\bar B_0)
t^\frac12\nu^{-\frac52}.
\end{equation}

\noindent Notice that
$$
\left\{x\in  D(A^{-1})\,,\,\alpha_\nu(x)< t\right\}\subset\mathcal
O_t,
$$
which yields, by \eqref{Gevrey_Eq2.3} and $\mu(H)=1$,
\begin{equation}\label{Gevrey_Eq2.2}
\mu\left(x\in  H\,,\,\alpha_\nu(x)\leq t\right)\leq  (K +1)(1+\bar
B_0) t^\frac12\nu^{-\frac52}.
\end{equation}
It is well-known that \eqref{Gevrey_Eq2.2} for any $t>0$ implies \eqref{Gevrey_Eq_Th_Gevrey_2},
 which yields Theorem \ref{Gevrey_Th_Gevrey}.

\subsection{Proof of Corollary \ref{Gevrey_Cor_Gevrey}}\label{Sec_Cor_Gevrey}

\

\noindent To deduce Corollary \ref{Gevrey_Cor_Gevrey} from Theorem \ref{Gevrey_Th_Gevrey}, it is sufficient to prove
that for any $(\alpha',\alpha,\beta',\beta)\in \left(0,\infty\right)^2\times\left(0,1\right]^2$ such that $\beta'<\beta$, we have
\begin{equation}\label{Gevrey_Eq4.1}
\norm{x}_{G(\alpha',\beta')}\leq
 \exp\left(c(\beta,\beta')\left(\alpha'\right)^\frac{\beta}{\beta-\beta'}
\left(\alpha\right)^{-\frac{\beta'}{\beta-\beta'}}\right)
\norm{x}_{G(\alpha,\beta)}.
\end{equation}
Indeed, \eqref{Gevrey_Eq4.1} implies that for any $\gamma\in\R^+_*$
$$
 \left(\ln^+ \norm{x}^2_{G(\alpha',\beta')}\right)^{\gamma}\leq
c_\gamma\left(c(\beta,\beta')+\left(\alpha'\right)^\frac{\gamma\beta}{\beta-\beta'}
\left(\nu\alpha_\nu(x)\right)^{-\frac{\gamma\beta'}{\beta-\beta'}}
+ \left(\ln^+ \norm{x}^2_{G(\nu \alpha_\nu(x),\beta)}\right)^{\gamma}\right),
$$
which yields  Corollary \ref{Gevrey_Cor_Gevrey} provided Theorem \ref{Gevrey_Th_Gevrey} is true.

\noindent We now establish \eqref{Gevrey_Eq4.1}. It follows from
arithmetic-geometric inequality that for any
 $k\in \Z^3$
\begin{equation}\label{Gevrey_Eq4.2}
\alpha' \abs{k}^{\beta'}\leq c(\beta,\beta')\left(\alpha'\right)^\frac{\beta}{\beta-\beta'}
\left(\alpha\right)^{-\frac{\beta'}{\beta-\beta'}}
+\alpha \abs{k}^\beta.
\end{equation}
 Recalling that
$$
\norm{x}_{G(\alpha',\beta')}^2=\sum_{k\in \Z^3}\abs{k}^2\exp\left(2\alpha'\abs{k}^{\beta'}\right)\abs{\hat x(k)}^2,
$$
  we infer \eqref{Gevrey_Eq4.1} from \eqref{Gevrey_Eq4.2}.


\end{document}